\documentclass[a4paper, 10 pt, conference]{article}                                                          
\newif\ifarticle
\articletrue

\usepackage{epsfig} 
\usepackage{booktabs}
\usepackage{amsmath} 
\usepackage{amssymb}  
\usepackage{enumerate}
\usepackage{url}

\newtheorem{thm}{Theorem}
\newtheorem{prp}{Proposition}

\newtheorem{ass}{Assumption}
\newtheorem{lem}{Lemma}

\newenvironment{pf}{\smallbreak\noindent{\it Proof. }}{\hfill$\Box$\smallbreak}

\usepackage{tikz}
\usepackage{pgfplots}
\usetikzlibrary{calc}
\newcommand{\marker}[2]
{
  \pgfgettransformentries{\myxscale}{\@tempa}{\@tempa}{\myyscale}{\@tempa}{\@tempa}
  \draw[thick] ($#1+0.08*(1/\myxscale,1/\myyscale)$)--($#1-0.08*(1/\myxscale,1/\myyscale)$);
  \draw[thick] ($#1+0.08*(-1/\myxscale,1/\myyscale)$)--($#1-0.08*(-1/\myxscale,1/\myyscale)$);
}

\newcommand{\id}{\rm{Id}}
\newcommand{\reals}{\mathbb{R}}

\newcommand{\alphanom}{\bar{\alpha}}
\newcommand{\alphak}{\alpha_k}
\newcommand{\fix}{{\rm{fix}}}

\newcommand{\cmt}[1]{}

\title{\LARGE \bf
Line Search for Generalized Alternating Projections
}

\author{Mattias F\"{a}lt and Pontus Giselsson}

\begin{document}
  \newlength\figureheight
  \newlength\figurewidth
  \newcommand{\figurecaptionreduction}{\vspace{-4mm}}
  
\maketitle
\thispagestyle{empty}
\pagestyle{empty}

\begin{abstract}
This paper is about line search for the generalized alternating projections (GAP) method.
This method is a generalization of the von Neumann alternating projections method, where
instead of performing alternating projections, relaxed projections are alternated.
The method can be interpreted as an averaged iteration of a nonexpansive mapping. Therefore,
a recently proposed line search method for such algorithms is applicable to GAP.
We evaluate this line search and show situations when the line search can be performed
with little additional cost. We also present a variation of the basic line search for GAP---the projected line search.
We prove its convergence and show that the line search condition
is convex in the step length parameter.
We show that almost all convex optimization problems can be solved using this approach
and numerical results show superior performance with both
the standard and the projected line search, sometimes by
several orders of magnitude, compared to the nominal method.
\end{abstract}

\newif\ifbackground
\backgroundtrue
\section{Introduction}
Alternating projections is a well known method for feasibilty problems,
where the objective is to find a point in the intersection of (convex) sets.
The method alternates projections onto the sets.
It was first introduced for half-spaces~\cite{GAP_Agmon} and later generalized to more general sets~\cite{BREGMAN1967200}.
In practice, the method is often quite slow.
A generalization to this method was proposed in~\cite{GAP_Gubin},
which is based on performing relaxed projections onto the sets instead of standard projections.
A relaxation parameter defines how far the relaxed projection should go towards or past the projection point.
Depending on the relaxation parameters,
it can be shown that the method is an averaged iteration of a nonexpansive mapping,
where fixed-points to the mapping correspond to solutions to the feasibility problem.
Many variations and extensions of this basic method has been proposed and studied,
with linear or sublinear convergence estimates~\cite{BauschkeFixPointsNonExpansive, BauschkeOnProjection}.
We present a framework for several of these generalizations and collect the relevant results.
The framework includes well known methods such as the alternating projections and
the generalized Douglas-Rachford algorithm for feasibility problems.
\cmt{Rewrite again? --- }These are first order methods that scale better with the number of variables compared to interior point algorithms
and usually have a low computational cost per iteration.
They are therefore useful for solving large-scale convex optimization problems,
where the computational complexity is too large for other algorithms.
The practical rate of convergence can however be slow and is dependent on preconditioning for good performance.
A good preconditioning can be hard to find and is usually problem specific\cmt{~Cite Pontus?}.

Line search is a well established concept in optimization
and is often used to improve practical performance of a method.
Typically, it assumes that a descent direction for the objective function is at hand,
and it accepts points with sufficient decrease and possibly some condition on the slope~\cite{Boyd2004, Nocedal}.
For averaged iterations of nonexpansive mappings,
descent directions are not obtained in general.
In the recent paper~\cite{gis_line_search},
a line search method that can be applied to averaged iterations was proposed.
The line search is performed in the direction of the fixed-point residual,
which is the direction obtained by applying the nonexpansive operator.\cmt{~Sounds weird?}
Instead of being based on objective function value decrease,
it relies on a decrease in the norm of the fixed-point residual.

The main contribution of this paper is an alternative to the basic line search for GAP---the projected line search.
This is developed for the case with two sets, where one is affine,
and we show that most convex optimization problems can be posed on this form.
The projected line search method performs line search, not in the residual direction,
but in its projection on the affine set.
We prove that the method converges,
and show that the line search condition is convex in the step length parameter.
We also present a numerical example that illustrates the properties of the methods,
and show that the projected line search can achieve superior performance.

\ifbackground
Section~\ref{sec:background} contains some background and notation.
\fi
In Section~\ref{sec:GAP} we present the generalized alternating projections algorithm and collect relevant results.
In Section~\ref{sec:ls} we show how the line search in~\cite{gis_line_search}
can be applied to this algorithm.
The projected line search is presented in Section~\ref{sec:projected_ls} together with some basic results.
An overview of how this method can be used to solve a large set of convex optimization problems is presented in Section~\ref{sec:cone},
and a numerical example is presented in Section~\ref{sec:numEx}.

\ifbackground
\section{Background and Notation}\label{sec:background}
The notation $\langle\cdot,\cdot\rangle$ is used for scalar product and
$\rm{Id}$ is the identity operator.
The fixed-points of an operator $T$ are denoted $\fix T$, i.e.~$\fix T=\{x\in\mathbb{R}^n\mid Tx=x\}$,
and the fixed-point residual $r(x)$ for a point $x$ is defined as $r(x) := Tx-x$.
An operator $T : \mathbb{R}^n\rightarrow \mathbb{R}^n$ is said to be nonexpansive if it satisfies $\|Tx-Ty\|_2\leq\|x-y\|_2$ for all $x,y\in\mathbb{R}^n$,
and an operator $S$ is $\alpha$-averaged, with $\alpha\in(0,1)$, if it can be written as $S=(1-\alpha)\rm{Id}+\alpha T$ for some nonexpansive $T$.
$\Pi_C$ is the orthogonal projection onto the closed, convex and nonempty set $C$, i.e. $\Pi_C(x)=\arg\min_{y\in C}(\|y-x\|_2)$.
\fi

\section{Generalized Alternating Projections}\label{sec:GAP}

Generalized alternating projections is an algorithm for finding a point in the intersection of $p$ sets $C_i$ with $i=1,\ldots,p$, i.e.,
to find a point $x\in C_1\cap\cdots\cap C_p$.
Thoughout this paper, we assume that that sets $C_i$ are nonempty, closed and convex, and that
\begin{align}
C_1\cap ...\cap C_p \neq \emptyset,
\label{eq:intersection_nonempty}
\end{align}
i.e., that a common feasible point exists.

To define the algorithm, we introduce the under ($\alpha\in(0,1)$) and over ($\alpha\in(1,2]$) relaxed projection on the set $C$ as follows:
\begin{align}
P_{C}^{\alpha}=(1-\alpha)\id+\alpha\Pi_C
\label{eq:rel_proj}
\end{align}
where $\alpha\in(0,2]$ and $\Pi_C$ is the orthogonal projection onto the set $C$. For $\alpha=1$,
we get the standard projection $P_C^1=\Pi_C$ and for $\alpha=2$,
we get the reflection $P_C^2=2\Pi_C-\id=:R_{C}$.
The relaxed projector is $\tfrac{\alpha}{2}$-averaged for $\alpha\in(0,2)$ and nonexpansive for $\alpha=2$,
since $\Pi_C$ is firmly nonexpansive, see~\cite[Corollary~4.29, Example~12.25, Proposition~12.27]{bauschkeCVXanal}.

The generalized alternating projections method (GAP) is:
\begin{align}
x^{k+1}=(1-\alpha)x^k+\alpha P_{C_{p}}^{\alpha_p}P_{C_{p-1}}^{\alpha_{p-1}}\cdots P_{C_1}^{\alpha_1}x^k.
\label{eq:gap}
\end{align}
For simplicity, we introduce the GAP operator $T$ as
\begin{align}
T = (1-\alpha)\id+\alpha P_{C_{p}}^{\alpha_p}P_{C_{p-1}}^{\alpha_{p-1}}\cdots P_{C_1}^{\alpha_1}
\label{eq:gap_op}
\end{align}
to arrive at the notationally more convenient iteration $x^{k+1}=Tx^k$ for~\eqref{eq:gap}.

The algorithm~\eqref{eq:gap} generalizes the classical alternating projections method, since
if $\alpha=\alpha_i=1$, we get
\begin{align*}
x^{k+1}=\Pi_{C_{p}}\Pi_{C_{p-1}}\cdots\Pi_{1}x^k.
\end{align*}

For $p=2$, the generalized Douglas-Rachford algorithm for feasibility problems~\cite{DouglasRachford, LionsMercier1979},
also falls under the formulation~\eqref{eq:gap} by letting $\alpha_1=\alpha_2=2$. Then
\begin{align*}
x^{k+1}=(1-\alpha)x^k+\alpha R_{C_2}R_{C_1}x^k
\end{align*}
where $R_C=2\Pi_C-\id$ is a reflector.
\begin{figure}
\xdef\thetainit{330}

\begin{center}
\centering
\begin{tikzpicture}


\begin{scope}[scale=7]
\clip ({cos(\thetainit)-0.35},{sin(\thetainit)-0.06}) rectangle (1.15,0.05);

\def\xline{1};
\draw[fill=gray!25!white] (0,0) circle (1);
\draw (\xline,-2)--(\xline,2);
\marker{(1,0)}{0.01}
\node[right] at (1,0) {$x^\star$};
\edef\xoneinit{cos(\thetainit)};
\edef\xtwoinit{sin(\thetainit)};
\edef\alphatwo{2};
\edef\alphanom{0.5};
\xdef\one{1}


\foreach \alphaone in {1,2} {

\edef\alphatwo{\alphaone};
\pgfmathparse{1/\alphaone};
\edef\alphanom{\pgfmathresult};

\pgfmathparse{1*\xoneinit}
\edef\xone{\pgfmathresult}
\pgfmathparse{1*\xtwoinit}
\edef\xtwo{\pgfmathresult}

\marker{(\xone,\xtwo)}{0.01};
\ifx\alphaone\one
\node[left] at (\xone,\xtwo) {$x^{1}$};
\fi
\edef\pxone{\xline};
\pgfmathparse{1*\xtwo}
\edef\pxtwo{\pgfmathresult}
\pgfmathparse{\alphatwo*\pxone+(1-\alphatwo)*\xone};
\edef\rxone{\pgfmathresult};
\pgfmathparse{\alphatwo*\pxtwo+(1-\alphatwo)*\xtwo};
\edef\rxtwo{\pgfmathresult};


\pgfmathparse{sqrt(\rxone*\rxone+\rxtwo*\rxtwo)};
\edef\rxnorm{\pgfmathresult};
\pgfmathparse{\rxone/\rxnorm};
\xdef\prxone{\pgfmathresult};
\pgfmathparse{\rxtwo/\rxnorm};
\xdef\prxtwo{\pgfmathresult};
\pgfmathparse{\alphaone*\prxone+(1-\alphaone)*\rxone};
\edef\rrxone{\pgfmathresult};
\pgfmathparse{\alphaone*\prxtwo+(1-\alphaone)*\rxtwo};
\edef\rrxtwo{\pgfmathresult};

\pgfmathparse{(1-\alphanom)*\xone+\alphanom*\rrxone};
\edef\nextxone{\pgfmathresult};
\pgfmathparse{(1-\alphanom)*\xtwo+\alphanom*\rrxtwo};
\edef\nextxtwo{\pgfmathresult};
\marker{(\nextxone,\nextxtwo)}{0.01};

\ifx\alphaone\one
\node[right] at (\nextxone,\nextxtwo) {$x_{{\rm{AP}}}^2$};
\else
\node[left] at (\nextxone,\nextxtwo) {$x_{{\rm{DR}}}^2$};
\fi

\draw[-latex,gray] (\xone,\xtwo)--(\rxone,\rxtwo);
\draw[-latex,gray] (\rxone,\rxtwo)--(\rrxone,\rrxtwo);
\draw[-latex,red] (\xone,\xtwo)--(\rrxone,\rrxtwo);

}

\end{scope}

\end{tikzpicture}

\end{center}
\caption{\label{fig:gap_ex}
Illustration of the generalized alternating projections for two different settings
on a 2-dimensional problem with intersection $x^\star$.
The first set is the vertical line, and the second is the shaded area.
The point $x^2_{\rm{AP}}$ is obtained by
an alternating projections step ($\alpha_1=\alpha_2=1$, $\alpha=1$)
and $x^2_{\rm{DR}}$ is obtained by a Douglas-Rachford step ($\alpha_1=\alpha_2=2$, $\alpha=0.5$).
The red arrows represent the residuals $P^{\alpha_2}_{C_2}P^{\alpha_1}_{C_1}x^1-x^1$
along which we will perform line search in Section~\ref{sec:ls}.
}
\end{figure}
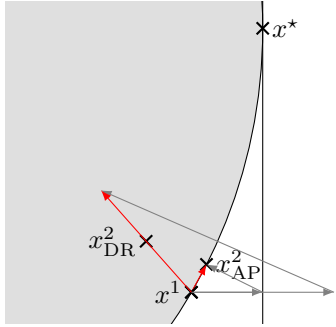
These two algorithms are illustrated for a simple 2-dimensional problem in Figure~\ref{fig:gap_ex}.

Below, we present some basic results on the algorithm~\eqref{eq:gap}.
Most of these are known but spread out in the literature
so we collect them here for convenience of the reader.
To this end, we let
\begin{align}
\beta :=\frac{\sum_{i=1}^p\tfrac{\alpha_i}{2-\alpha_i}}{1+\sum_{i=1}^p\tfrac{\alpha_i}{2-\alpha_i}},
\label{eq:beta}
\end{align}
and state the following assumptions on $\alpha$.
\begin{ass}
Suppose that either of the following holds:
\begin{enumerate}[{\it{A}}1.]
\item $\alpha\in(0,\tfrac{1}{\beta})$ with $\beta$ in~\eqref{eq:beta} and that $\alpha_i\in(0,2)$ for $i=1,\ldots,p$
\item $\alpha\in(0,1)$ and $\alpha_i\in(0,2]$ for $i=1,\ldots,p$ with at most one $\alpha_i=2$
\item $\alpha\in(0,1)$ and $p=2$ with $\alpha_1=\alpha_2=2$.
\end{enumerate}
\label{ass:alpha}
\end{ass}

These assumptions imply that the GAP operator $T$ is averaged. This is shown next.
\begin{prp}\label{prop:prop1}
Suppose that Assumption~\ref{ass:alpha} with case A1 holds.
Then the GAP operator $T$ in \eqref{eq:gap_op} is averaged with constant $\alpha\beta\in(0,1)$,
with $\beta$ in~\eqref{eq:beta}. Suppose that Assumption~\ref{ass:alpha} with case A2 or A3 holds.
Then $T$ is averaged with constant $\alpha\in(0,1)$.
\end{prp}
A proof is found in the Appendix.

Next, we show a result on the fixed-point set of the GAP operator in~\eqref{eq:gap_op}.
\begin{prp}\label{prp:fixT}
Suppose that Assumption~\ref{ass:alpha} holds with case A1 or A2 and that $C_1\cap\cdots\cap C_p\neq\emptyset$.
Then $\fix T = C_1\cap\cdots\cap C_p$,
where $T$ is the operator in~\eqref{eq:gap_op}.
\end{prp}
A proof is found in the Appendix.

The main convergence result for the algorithm now follows directly from~\cite[Theorem 5.14]{bauschkeCVXanal}
under assumption A1 or A2 since $T$ is averaged
and its fixed-point set is $C_1\cap\cdots\cap C_p$.
\begin{prp}
Suppose that Assumption~\ref{ass:alpha} holds  with case A1 or A2 and that $C_1\cap\cdots\cap C_p\neq\emptyset$.
The fixed-point residuals $r(x^k)$  converge to $0$
and the iterates $x^k$ converge to a point in the intersection $C_1\cap\cdots\cap C_p$,
as $k\rightarrow \infty$ in algorithm~\eqref{eq:gap}.
\end{prp}

Algorithm~\eqref{eq:gap} with case A3 in Assumption~\ref{ass:alpha}
corresponds to generalized Douglas-Rachford applied to feasibility problems.
The properties in this case are slightly different but well known,
and we summarize them below~\cite[Proposition 25.1, Theorem 25.6]{bauschkeCVXanal}.
\begin{prp}
Suppose that Assumption~\ref{ass:alpha} holds  with case A3 and that $C_1\cap C_2\neq\emptyset$,
then the fixed-point set satisfies $\Pi_{C_1}\fix T=C_1\cap C_2$.
Additionally, the fixed-point residuals $r(x^k)$ in algorithm~\eqref{eq:gap} converge to $0$ as $k\rightarrow \infty$ and
the iterates $x^k$ converge to a point $x$
such that $\Pi_{C_1}x\in C_1\cap C_2$.
\end{prp}

We see that we need to monitor the sequence $\Pi_{C_1}x^k$ to find a feasible point in the Douglas-Rachford case.
For other choices of $\alpha_i$, it is also typically better to monitor
(one of) the projected sequences than $x^k$ in
\eqref{eq:gap} to faster find an intersection point (up to numerical accuracy).

\section{Line search}\label{sec:ls}
A method for applying line search on algorithms based on iterating averaged operators was recently proposed in~\cite{gis_line_search}.
The method was shown to often improve practical convergence.
In this section we describe how the method can be applied to generalized alternating projections.
We also repeat the result that under some assumptions on the sets,
the line search can be carried out with little additional cost compared to a basic iteration.

The line search algorithm can be applied to averaged iterations of the form
\begin{align}
x^{k+1}=(1-\alpha)x^k+\alpha Sx^k = x^k + \alpha (Sx^k-x^k),
\label{eq:averaged_iter} \end{align}
where $\alpha\in(0,1)$ and $S~:~\reals^n\to\reals^n$ is nonexpansive.
GAP is precisely on this form with $S=P_{C_{p}}^{\alpha_p}P_{C_{p-1}}^{\alpha_{p-1}}\cdots P_{C_1}^{\alpha_1}$.
The second expression in~\eqref{eq:averaged_iter}
shows that an averaged iteration performs a step with length $\alpha$ in the residual direction $r(x)=Sx-x$.
We call this the nominal step $\bar x^k := x^k+\alpha r(x^k)$.
The residual direction is illustrated in Figure~\ref{fig:gap_ex}.

The line search scheme presented in~\cite{gis_line_search}, suggests to perform line search in the residual direction.
To do this, the $\alpha$ that multiplies the residual direction should be chosen on-line.
The algorithm with line search can be written as:
\begin{align}
x^{k+1} &:= x^k + \alphak (Sx^k-x^k) := x^k + \alphak r(x^k)
\end{align}
where the line search parameter $\alphak$ must satisfy either $\alphak=\alpha$, i.e.,
we take the nominal step $\bar x^k$, or
$\alphak\in (\alpha, \alpha^{\max}]$ is such that
\begin{align}
\|r(x^{k+1})\|_2 \leq(1-\epsilon)\|r(\bar x^k)\|_2
\label{eq:ls_test}
\end{align}
where $\epsilon\in(0,1)$ and $\alpha^{\max} \geq \alpha$ are fixed algorithm parameters.
To accept a step length $\alphak$ in the line search,
the residual $r(x)$ should be smaller for the next iterate $x^{k+1}$
than for the nominal step $\bar x^k$.
This preserves the non-increasing property of the fixed-point residual $\|r(x^{k+1})\|$,
even when line search is used.
As shown in~\cite{gis_line_search}, this is enough to, e.g., guarantee convergence of the residual sequence.
An appropriate $\alphak$ can for example be selected using a simple forward or backward tracking.

The following form of the algorithm shows which computations are needed in each iteration:
\begin{align}
\label{eq:rk}r^k&:=Sx^k-x^k\\
\label{eq:xknom}\bar x^k&:=x^k+\alpha r^k\\
\label{eq:rknom}\bar r^k&:=S \bar x^k - \bar x^k\\
\label{eq:xk+1}x^{k+1} &:=x^k+\alphak r^k,
\end{align}
where $S=P_{C_{p}}^{\alpha_p}P_{C_{p-1}}^{\alpha_{p-1}}\cdots P_{C_1}^{\alpha_1}$.
The criterion for line search, i.e.~accepting $\alphak\neq\alpha$ in~\eqref{eq:xk+1}, can be written
\begin{align}
\|r^{k+1}\|_2 = \|Sx^{k+1}-x^{k+1}\|_2 \leq(1-\epsilon)\|\bar r^k\|_2
\label{eq:ls_test1}
\end{align}
where $x^{k+1} = x^k+\alphak r^k$ as seen in~\eqref{eq:xk+1}.
This general form of the algorithm reveals that we need to compute $S(x^k+\alphak r^k)$
for each candidate $\alphak$ to verify~\eqref{eq:ls_test},
as well as calculating $S\bar x^k$ each iteration.
So, to evaluate a candidate point in the line search is roughly as costly as performing one basic step in the algorithm.
This may or may not be too costly compared to what is saved due to the line search.

In the following,
we will show that sometimes many candidate points can be evaluated
in the line search with little additional cost.
In the case where the sets $C_n,\dots,C_1$ are affine, i.e. $C_i = \{x\in\reals~|~A_ix=b_i\}$,
the projection $z=\Pi_{C_i}x$ is affine and given by the solution to the KKT conditions of the projection:
\begin{align*}
  \begin{bmatrix}
    I & A_i^T\\
    A_i & 0
  \end{bmatrix}
  \begin{bmatrix} z\\ \lambda \end{bmatrix}
  = \begin{bmatrix} x\\ b_i \end{bmatrix}.
\end{align*}
The relaxed projections will therefore also be affine:
\begin{align*}
P^{\alpha_i}_{C_i}x&=(1-\alpha_i)x+\alpha_i \begin{bmatrix}I&0\end{bmatrix}\begin{bmatrix} I & A_i^T\\
A_i & 0\end{bmatrix}^{-1}\begin{bmatrix}x\\b_i\end{bmatrix}.
\end{align*}
It follows that the composition $P^{\alpha_n}_{C_n}\cdots P^{\alpha_1}_{C_1}$ is affine,
and GAP~\eqref{eq:gap} can be written as:
\begin{align*}
x^{k+1}=(1-\alpha)x^k+\alpha S_2S_1
\end{align*}
where $S_1x=P^{\alpha_n}_{C_n}\cdots P^{\alpha_1}_{C_1}x=: Fx+h$,
with $F$ and $h$ implicitly defined, and $S_2 = P^{\alpha_p}_{C_p}\cdots P^{\alpha_{n+1}}_{C_{n+1}}$.
The following iterations show that several candidate $\alpha_k$ can be tested,
without multiple evaluations of $S_1$~\cite{gis_line_search}:
\begin{align}
\label{eq:rk_comp2}r^k&:=S_2(Fx^k+h)-x^k\\
\label{eq:xknom_comp2}\bar x^k&:=x^k+\alpha r^k\\
\label{eq:rknom_comp2}\bar r^k&:=S_2\left(Fx^k +h +\alpha Fr^k\right) - \bar x^k\\
\label{eq:xk+1_comp} x^{k+1}&:=x^k+\alphak r^k
\end{align}
where $\alphak$ is selected so that
\begin{align}
\|S_2(Fx^{k+1}+h)-x^{k+1}\|_2 \leq(1-\epsilon)\|\bar r^k\|_2.
\label{eq:ls_test_comp}
\end{align}
The computed quantity $Fx^{k+1}=Fx^{k}+\alphak Fr^{k}$
is then reused in~\eqref{eq:rk_comp2},~\eqref{eq:rknom_comp2} and~\eqref{eq:ls_test_comp} in the following iteration.
Therefore, we only need to compute $Fx^0$ and $Fr^k$ for all $k$ to evaluate any number of candidate
$\alphak$ in any number of line searches.
If the cost of applying $S_2$ is negligible,
then the line search will result in no significant increase in computation per iteration.

\section{Projected line search}\label{sec:projected_ls}
In this section we present an alternative to the standard line search, that we call projected line search.
We present this line search for feasibility problems with two sets, $C_1=C$ and $C_2=D$, where C is affine.
This method does not select the next iterate in the direction of the residual but rather along its projection on the affine set.

The proposed algorithm, with $S=P^{\alpha_2}_DP^{\alpha_1}_C$, is:
\begin{align}
\label{eq:rkAlt}r^k&:=Sx^k-x^k\\
\label{eq:xknomAlt}\bar x^k&:=x^k+\alpha r^k\\
\label{eq:rknomAlt}\bar r^k&:=S \bar x^k - \bar x^k
\end{align}
where the next step is to either take a nominal step:
\begin{subequations}
\begin{align}
  \label{eq:xk+1Alt1}\quad x^{k+1} &:=x^k+\alpha r^k=\bar x^k\\
  \intertext{or line search is performed:}
  \quad\label{eq:xk+1Alt2}x^{k+1} &:=\Pi_C(x^k+\alphak r^k).
\end{align}
\end{subequations}
To accept the line search in~\eqref{eq:xk+1Alt2},
the line search parameter $\alphak\in (\alpha, \alpha^{\max}]$ must satisfy the following constraint,
where $i_{LS}$ is the index when the last line search was performed,
and $r(x^{i_{\rm{LS}+1}})= Sx^{i_{\rm{LS}+1}}-x^{i_{\rm{LS}+1}}$ is the residual at the following iteration:
\begin{align}
\|r(x^{k+1})\|_2 \leq(1-\epsilon)\|r(x^{i_{\rm{LS}+1}})\|_2.
\label{eq:ls_test2}
\end{align}

Compared to the algorithm with basic line search method,
the difference is that the candidate points in the projected line search are projected onto the set $C$.
The test for accepting a line search is also different.
Instead of comparing the norm of the next residual $r(x^{k+1})$ to the residual of the nominal step $r(\bar x^k)$,
we compare it to the residual in the last step that was chosen by a line search, $r(x^{i_{\rm{LS}+1}})$.
The reason for comparing the residual to the iterate $x^{i_{\rm{LS}+1}}$
is that the projected line search often increases the residual compared to the nominal step $\bar x^k$.
However, by ensuring that the residual $r(x^{k+1})$ is smaller than $r(x^{i_{\rm{LS}+1}})$,
we can guarantee that it will eventually decrease.
This is proven for general line search schemes in~\cite{gis_line_search} and we state it for the projected line search below.
\begin{thm}\label{thm:proj-converge}
Assume that Assumption~\ref{ass:alpha} holds and the projected line search algorithm~\eqref{eq:rkAlt}-\eqref{eq:xk+1Alt2}
is used with line search criteria~\eqref{eq:ls_test2}.
Then the fixed-point residual $r(x^k)=Sx^k-x^k$ will converge to $0$ as $k\rightarrow \infty$.
\end{thm}
We now show two additional properties of the projected line search.
\begin{thm}\label{thm:proj-cheap}
Assume that the set $C$ is affine,
then the projected line search condition~\eqref{eq:ls_test2} is convex in the step length $\alphak$
and the norm of the residual simplifies to $\|r(x^{k+1})\|_2= \alpha_2{\rm{dist}}_D(x^{k+1})$.
\end{thm}
A proof is found in the Appendix.

This result implies that we are not restricted to forward or back tracking schemes when finding step length $\alphak$.
We can perform, e.g., golden section search on the line search condition, or bisection on its gradient.
This way the number of candidate points to be evaluated, before a good/optimal point is found, can be reduced.
The theorem also illustrates that minimizing the left hand side of the line search condition is equivalent to
minimizing the distance between the two sets along the line search direction.
This gives an intuitive explanation to why this is a reasonable objective function.

We showed that the standard line search could be performed without significant extra computational cost in some cases,
the same is true for the projected line search.
If $C$ is affine, then $x^{k+1}=\Pi_C(x^k+\alphak r^k)$ can be evaluated for several $\alphak$
without any significant cost since the linear parts of $\Pi_C x^k$ and $\Pi_C r^k$ are known from previous steps,
in the same way as for the basic line search.
From Theorem~\ref{thm:proj-cheap} we know that evaluating the residual simplifies to evaluating the distance from $x^{k+1}$ the set $D$.
Thus, if $C$ is affine and $D$ is relatively cheap to project on, the line search will incur no significant cost.

\section{Cone programming}\label{sec:cone}

Many convex optimization problems, including LP, QP, SOCP, SDP,
and in particular problems that can be solved using optimization modeling interfaces such as CVX~\cite{cvx_v3},
CVXPY~\cite{cvxpy}, Convex.jl~\cite{convexjl},
can be written as cone programs of the form
\begin{align*}
\begin{tabular}{ll}
minimize & $c^Tx$\\
subject to & $Ax+s=b$\\
&$s\in\mathcal{K}$
\end{tabular}
\end{align*}
where $\mathcal{K}$ is a (product of) nonempty, closed and convex cones. The dual of this problem is
\begin{align*}
\begin{tabular}{ll}
maximize & $-b^Ty$\\
subject to & $-A^Ty=c$\\
&$y\in\mathcal{K}^*$.
\end{tabular}
\end{align*}
Assuming strong duality, we get that $c^Tx+b^Ty=0$, so the dual and primal problems can be embedded into the following primal-dual feasibility problem
\begin{align*}
\begin{tabular}{ll}
find & $(x,s,y)$\\
subject to & $\begin{bmatrix} A & I & 0\\ 0 & 0 & -A^T\\ c^T & 0&b^T \end{bmatrix}\begin{bmatrix}x\\s\\y\end{bmatrix}=\begin{bmatrix}b\\c\\0\end{bmatrix}$ \\
& $(s,y)\in\mathcal{K}\times\mathcal{K}^*$.
\end{tabular}
\end{align*}
This is a feasibility problem with one affine subspace and one product of convex cones.
There are many other ways to construct a feasibility problem with an affine subspace and a product of convex cones.
One example is the homogeneous self-dual embedding which is often used in interior-point methods~\cite{ye1994}
and in the first-order optimization solver SCS~\cite{SCS}.
Therefore, most convex optimization problems (at least those that can be posed as cone programs)
can be solved using GAP,
with one affine subspace and one product of convex cones.
This is precisely the formulation for which the basic line search
and the projected line search can be carried out with little additional cost
and where the line search condition for the projected line search is convex in the line search parameter.

\section{Numerical example}\label{sec:numEx}

In this section we demonstrate the performance improvements of GAP when line search is used.
We consider the following problem
\begin{align}
\begin{tabular}{ll}
find & $z$\\
such that &$Q(z-p)=0 $\\
&$z\geq 0$,
\end{tabular}
\label{eq:example-positive}
\end{align}
where $p=10^{-7}\mathbf{1}$ to guarantee feasibility of the problem,
and $Q\in\reals^{50\times 100}$ is randomly generated with independent normally distributed elements with unit variance and zero mean.

We define two sets as $C=\{z\ |\ Q(z-p)=0 \}$ and $D=\{z\ |\ z\geq 0 \}$.
Depending on $Q$, the feasible set $C\cap D$ may be very small or consist of infinitely long rays in the nonnegative orthant.
For this particular problem, $Q$ is generated such that no ray in the affine set lies completely in the nonnegative orthant.
Therefore, the intersection is relatively small.

As a termination criteria we use the following high accuracy requirement:
\begin{align*}
  \|Q(z-p)\|_2 & \leq 10^{-10}\\
  z & \geq 0,
\end{align*}
and we let $\alpha_1=\alpha_2\in[1,2]$ and $\alpha=0.85/\beta$, with $\beta$ in~\eqref{eq:beta}.

As proposed in~\cite{gis_line_search}, we do not perform line search in each iteration, but use the rule
\begin{align}
  \label{eq:trigger}\frac{\langle r^{k}, \bar r^{k}\rangle}{\|r^{k}\|_2\|\bar r^{k}\|_2}<1-10^{-4}
\end{align}
to trigger it.
The reason is that this often improves performance more than if line search is used in every iteration.
The rationale behind the rule is that a large $\alphak$ can often be accepted when the iterates are moving along a straight line,
i.e.~when the angle between consecutive iterates is small.
Numerical experiments suggest that consecutive iterates along a line seems to coincide with slow convergence,
further motivating the use of line search when this occurs.

Both the basic and the projected line search is performed using a simple forward-tracking with a factor $1.4$,
and the result for different $\alpha_1=\alpha_2$
are shown in Figure~\ref{fig:gap_exiters}.
The norm of the residual for each iteration is shown in Figure~\ref{fig:gap_exnormk}, for two different $\alpha_1=\alpha_2$.

\begin{figure}
{
\setlength{\figurewidth}{0.86\linewidth}
\setlength{\figureheight }{5cm}
{\scriptsize

\begin{tikzpicture}

\definecolor{color1}{rgb}{0.888873524661399,0.43564916971766,0.278122978814386}
\definecolor{color0}{rgb}{0,0.605603136637209,0.978680130252957}
\definecolor{color2}{rgb}{0.242224297852199,0.643275093157631,0.304448651534115}

\begin{axis}[
xlabel={$ a_1=a_2 $},
ylabel={Iterations},
xmin=1, xmax=2,
ymin=40, ymax=10000000,
ymode=log,
scale only axis,
width=\figurewidth,
height=\figureheight,
ytick={1,10,100,1000,10000,100000,1000000,10000000,100000000},
y label style={at={(axis description cs:0.05,.5)}},
xmajorgrids,
ymajorgrids,
legend entries={{No LS},{Standard LS},{Projected LS}},
legend cell align={left}
]
\addplot [color0]
table {%
1 5386
1.01 5316
1.02 5247
1.03 5177
1.04 5109
1.05 5041
1.06 4973
1.07 4905
1.08 4838
1.09 4770
1.1 4704
1.11 4638
1.12 4572
1.13 4507
1.14 4442
1.15 4377
1.16 4313
1.17 4249
1.18 4185
1.19 4122
1.2 4059
1.21 3996
1.22 3934
1.23 3872
1.24 3811
1.25 3750
1.26 3688
1.27 3627
1.28 3567
1.29 3507
1.3 3447
1.31 3388
1.32 3329
1.33 3271
1.34 3213
1.35 3156
1.36 3098
1.37 3042
1.38 2985
1.39 2929
1.4 2873
1.41 2817
1.42 2761
1.43 2706
1.44 2651
1.45 2597
1.46 2542
1.47 2489
1.48 2435
1.49 2382
1.5 2328
1.51 2276
1.52 2223
1.53 2171
1.54 2119
1.55 2067
1.56 2016
1.57 1964
1.58 1913
1.59 1862
1.6 1811
1.61 1761
1.62 1711
1.63 1660
1.64 1610
1.65 1560
1.66 1510
1.67 1460
1.68 1411
1.69 1361
1.7 1311
1.71 1262
1.72 1213
1.73 1163
1.74 1114
1.75 1064
1.76 1015
1.77 965
1.78 916
1.79 866
1.8 816
1.81 765
1.82 715
1.83 663
1.84 611
1.85 557
1.86 501
1.87 440
1.88 372
1.89 252
1.9 113
1.91 120
1.92 129
1.93 143
1.94 159
1.95 185
1.96 222
1.97 282
1.98 396
1.99 727
2 8423438
};
\addplot [color1]
table {%
1 889
1.01 895
1.02 974
1.03 883
1.04 830
1.05 956
1.06 912
1.07 854
1.08 824
1.09 868
1.1 960
1.11 821
1.12 899
1.13 945
1.14 822
1.15 838
1.16 788
1.17 801
1.18 807
1.19 738
1.2 804
1.21 862
1.22 826
1.23 656
1.24 800
1.25 857
1.26 727
1.27 764
1.28 794
1.29 744
1.3 794
1.31 715
1.32 702
1.33 701
1.34 701
1.35 669
1.36 663
1.37 683
1.38 696
1.39 631
1.4 618
1.41 624
1.42 649
1.43 606
1.44 583
1.45 616
1.46 540
1.47 544
1.48 533
1.49 531
1.5 582
1.51 558
1.52 542
1.53 519
1.54 519
1.55 470
1.56 454
1.57 465
1.58 442
1.59 449
1.6 422
1.61 405
1.62 404
1.63 423
1.64 373
1.65 396
1.66 365
1.67 403
1.68 384
1.69 390
1.7 387
1.71 365
1.72 343
1.73 342
1.74 346
1.75 306
1.76 324
1.77 319
1.78 301
1.79 317
1.8 307
1.81 261
1.82 272
1.83 287
1.84 288
1.85 281
1.86 256
1.87 259
1.88 242
1.89 218
1.9 102
1.91 105
1.92 98
1.93 110
1.94 121
1.95 114
1.96 93
1.97 154
1.98 150
1.99 235
2 8423390
};
\addplot [color2]
table {%
1 857
1.01 986
1.02 816
1.03 925
1.04 888
1.05 832
1.06 857
1.07 706
1.08 841
1.09 800
1.1 846
1.11 831
1.12 734
1.13 648
1.14 736
1.15 838
1.16 784
1.17 759
1.18 647
1.19 623
1.2 656
1.21 641
1.22 567
1.23 486
1.24 640
1.25 534
1.26 509
1.27 633
1.28 655
1.29 527
1.3 677
1.31 459
1.32 478
1.33 587
1.34 470
1.35 429
1.36 564
1.37 350
1.38 547
1.39 572
1.4 442
1.41 474
1.42 524
1.43 589
1.44 472
1.45 388
1.46 365
1.47 447
1.48 326
1.49 403
1.5 385
1.51 350
1.52 441
1.53 388
1.54 342
1.55 360
1.56 392
1.57 404
1.58 354
1.59 426
1.6 392
1.61 308
1.62 316
1.63 280
1.64 300
1.65 368
1.66 260
1.67 381
1.68 273
1.69 298
1.7 326
1.71 230
1.72 223
1.73 253
1.74 244
1.75 234
1.76 242
1.77 274
1.78 193
1.79 191
1.8 220
1.81 197
1.82 188
1.83 169
1.84 171
1.85 211
1.86 183
1.87 177
1.88 207
1.89 134
1.9 76
1.91 57
1.92 69
1.93 69
1.94 65
1.95 52
1.96 68
1.97 54
1.98 70
1.99 70
2 91
};

\end{axis}
\end{tikzpicture}}
}
\figurecaptionreduction
\caption{\label{fig:gap_exiters}Number of iterations to solve problem~\eqref{eq:example-positive} for different $\alpha_1,\alpha_2$,
with and without line search.}
\end{figure}
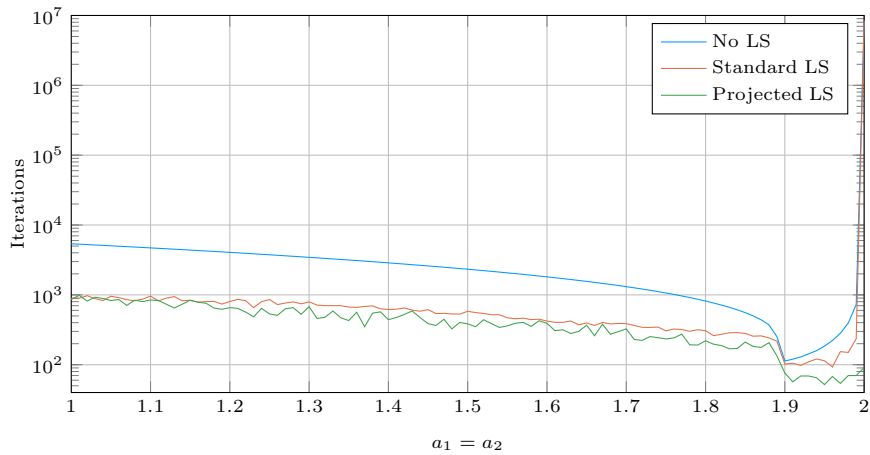
\begin{figure}
{
\setlength{\figurewidth}{\linewidth}
\setlength{\figureheight }{4cm}
{\scriptsize\input{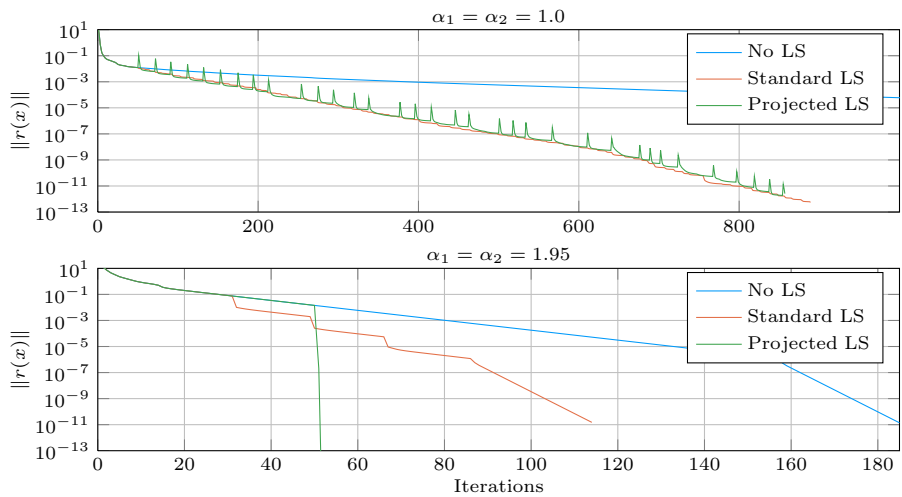}}
}
\caption{\label{fig:gap_exnormk} Norm of the residual for each iteration when solving
problem~\eqref{eq:example-positive} with different settings.
It can be noted that the norm is strictly decreasing both without line search and with the standard line search.
The peaks for the projected line search represents when a candidate $\alphak$ was accepted,
which sometimes result in a temporary increase in the norm
due to the constructed line search condition.}
\end{figure}
\bibliographystyle{plain}

\cmt{Maybe reformulate}Without line search, it is clear from Figure~\ref{fig:gap_exiters} that the choice $\alpha_1=\alpha_2=1$ and $\alpha_1=\alpha_2=2$,
corresponding to the Alternating Projections and Douglas-Rachford splitting respectively,
are far from the optimal choice.
They require approximately $5000$ and $8\cdot 10^6$ iterations
compared to only $113$ iterations for the optimal $\alpha_1=\alpha_2$.
However, this behavior does not apply to all problems,
for example, in many cases the number of iterations seems to be monotonically decreasing with larger $\alpha_1=\alpha_2$.

Figure~\ref{fig:gap_exiters} also reveals that the basic line search can considerably improve performance,
especially for low choices of $\alpha_1,\alpha_2$,
while the improvement for larger values is more modest.
However, the projected line search performs considerably better,
even for large $\alpha_1=\alpha_2$,
with only $52$ iterations for the optimal $\alpha_1=\alpha_2$.
In particular,
it decreases the iterations for $\alpha_1=\alpha_2=2$ with more than a factor $10^5$.

\begin{figure}
  \centering
  \begin{tikzpicture}[every node/.style={scale=0.65}]

\definecolor{color1}{rgb}{0.3,0.8,0.3}
\definecolor{color0}{rgb}{0.8,0.3,0.3}
\definecolor{color2}{rgb}{0.3,0.3,0.8}

\begin{axis}[
xmin=0.50, xmax=3.5,
ymin=0, ymax=4.2,
xtick={1,2,3},
xticklabels={No LS,Basic LS,Projected LS},
xticklabel style={anchor=north, yshift=-0.6em},
ytick={1,2,3,4},
yticklabels={$10^1$,$10^2$,$10^3$,$10^4$},
ymajorgrids,
legend entries={{} {Iterations},{} {LS Triggered},{} {LS Accepted}},
legend cell align={left},
ybar,
clip=false
]
\addplot[fill=color0] coordinates {(1,0)};
\addplot[fill=color1] coordinates {(1,0)};
\addplot[fill=color2] coordinates {(1,0)};
\path [draw=black, fill=color0] (axis cs:0.75,3.73126634907549)
--(axis cs:0.75,0)
--(axis cs:0.95,0)
--(axis cs:0.95,3.73126634907549)
--(axis cs:0.75,3.73126634907549);
\node at (axis cs:0.85,3.73126634907549) [anchor=south] {$5386$};

\path [draw=black, fill=color0] (axis cs:1.75,2.94890176097021)
--(axis cs:1.75,0)
--(axis cs:1.95,0)
--(axis cs:1.95,2.94890176097021)
--(axis cs:1.75,2.94890176097021);
\node at (axis cs:1.85,2.94890176097021) [anchor=south] {$889$};

\path [draw=black, fill=color0] (axis cs:2.75,2.9329808219232)
--(axis cs:2.75,0)
--(axis cs:2.95,0)
--(axis cs:2.95,2.9329808219232)
--(axis cs:2.75,2.9329808219232);
\node at (axis cs:2.85,2.9329808219232) [anchor=south] {$857$};

\path [draw=black, fill=color0] (axis cs:1.05,2.26717172840301)
--(axis cs:1.05,0)
--(axis cs:1.25,0)
--(axis cs:1.25,2.26717172840301)
--(axis cs:1.05,2.26717172840301);
\node at (axis cs:1.15,2.26717172840301) [anchor=south] {$185$};

\path [draw=black, fill=color0] (axis cs:2.05,2.05690485133647)
--(axis cs:2.05,0)
--(axis cs:2.25,0)
--(axis cs:2.25,2.05690485133647)
--(axis cs:2.05,2.05690485133647);
\node at (axis cs:2.15,2.05690485133647) [anchor=south] {$114$};

\path [draw=black, fill=color0] (axis cs:3.05,1.7160033436348)
--(axis cs:3.05,0)
--(axis cs:3.25,0)
--(axis cs:3.25,1.7160033436348)
--(axis cs:3.05,1.7160033436348);
\node at (axis cs:3.15,1.7160033436348) [anchor=south] {$52$};


\path [draw=black, fill=color1] (axis cs:1.70,2.21748394421391)
--(axis cs:1.7,0)
--(axis cs:1.9,0)
--(axis cs:1.9,2.21748394421391)
--(axis cs:1.7,2.21748394421391);
\node at (axis cs:1.7,2.21748394421391) [anchor=south east] {$165$};

\path [draw=black, fill=color1] (axis cs:2.7,1.53147891704226)
--(axis cs:2.7,0)
--(axis cs:2.9,0)
--(axis cs:2.9,1.53147891704226)
--(axis cs:2.7,1.53147891704226);
\node at (axis cs:2.7,1.53147891704226) [anchor=south east] {$34$};


\path [draw=black, fill=color1] (axis cs:2.1,0.778151250383644)
--(axis cs:2.1,0)
--(axis cs:2.3,0)
--(axis cs:2.3,0.778151250383644)
--(axis cs:2.1,0.778151250383644);
\node at (axis cs:2.3,0.778151250383644) [anchor=west] {$6$};

\path [draw=black, fill=color1] (axis cs:3.1,1.32221929473392)
--(axis cs:3.1,0)
--(axis cs:3.3,0)
--(axis cs:3.3,1.32221929473392)
--(axis cs:3.1,1.32221929473392);
\node at (axis cs:3.3,1.32221929473392) [anchor=west] {$21$};


\path [draw=black, fill=color2] (axis cs:1.65,2.12710479836481)
--(axis cs:1.65,0)
--(axis cs:1.85,0)
--(axis cs:1.85,2.12710479836481)
--(axis cs:1.65,2.12710479836481);
\node at (axis cs:1.65,2.12710479836481) [anchor=east] {$134$};

\path [draw=black, fill=color2] (axis cs:2.65,1.53147891704226)
--(axis cs:2.65,0)
--(axis cs:2.85,0)
--(axis cs:2.85,1.53147891704226)
--(axis cs:2.65,1.53147891704226);
\node at (axis cs:2.65,1.53147891704226) [anchor=north east] {$34$};


\path [draw=black, fill=color2] (axis cs:2.15,0.477121254719662)
--(axis cs:2.15,0)
--(axis cs:2.35,0)
--(axis cs:2.35,0.477121254719662)
--(axis cs:2.15,0.477121254719662);
\node at (axis cs:2.35,0.477121254719662) [anchor=west] {$3$};

\path [draw=black, fill=color2] (axis cs:3.15,0.301029995663981)
--(axis cs:3.15,0)
--(axis cs:3.35,0)
--(axis cs:3.35,0.301029995663981)
--(axis cs:3.15,0.301029995663981);
\node at (axis cs:3.35,0.301029995663981) [anchor=west] {$2$};


%
%

%
%

\node at (axis cs:1,0) [anchor=north east] {$\alpha=1$};
\node at (axis cs:1,0) [anchor=north west] {$\alpha=1.95$};

\node at (axis cs:2,0) [anchor=north east] {$\alpha=1$};
\node at (axis cs:2,0) [anchor=north west] {$\alpha=1.95$};

\node at (axis cs:3,0) [anchor=north east] {$\alpha=1$};
\node at (axis cs:3,0) [anchor=north west] {$\alpha=1.95$};

\end{axis}

\end{tikzpicture}
  \caption{Number of times the line search was triggered and accepted for different algorithms and settings.}
  \label{tab:LS}
\end{figure}
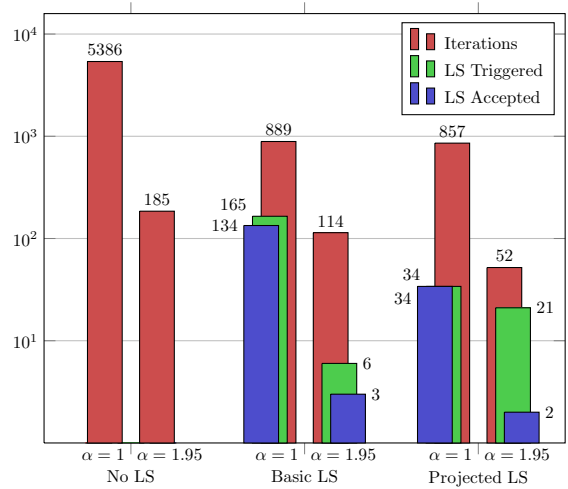

So far, we only compared the number of iterations for the different methods.
Since the line search methods have a higher per iteration cost,
we now evaluate what is actually gained by performing line search.
We focus on the two cases with $\alpha_1=\alpha_2=1.0$ and $1.95$. Figure~\ref{tab:LS}
shows the number of times the trigger criterion for line search~\eqref{eq:trigger}
was satisfied for the standard and projected line search.
It also shows how many times the line search found a point that satisfied the
corresponding criterion for line search acceptance,~i.e.~\eqref{eq:ls_test} and~\eqref{eq:ls_test2}.
The number of evaluated candidate points (i.e.~different $\alphak$)
averaged around $10$ for each line search attempt, with a maximum of $18$.
Since $C=\{z\ |\ Q(z-p)=0 \}$ is affine, only one extra projection on $C$ was needed
for line search (to initialize the algorithm).
To evaluate the acceptance criterion,~\eqref{eq:ls_test} or~\eqref{eq:ls_test2}, a few vector operations and
one projection onto $D=\{z\ |\ z\geq 0 \}$ is needed for each $\alphak$.
But projecting onto $D$ is simply a $\max$-operation and is thus very cheap. \cmt{maybe compare the cost of $\Pi_C$ and $\Pi_D$?}

\section{Conclusions}
We have shown that a recently proposed line search~\cite{gis_line_search}
is applicable to the generalized alternating projections method.
We have also proposed an alternative line search method for GAP,
the projected line search.
Furthermore, we have shown that the line search condition for the projected line search is convex in the step length parameter.
Both line search methods were evaluated on a feasibility problem,
and showed significant performance improvements compared to the nominal method.

\section*{Acknowledgment}
Both authors are financially supported by the Swedish Foundation for Strategic
Research and members of the LCCC Linneaus Center at Lund University.

\bibliography{references.bib}

\begin{thebibliography}{10}

\bibitem{GAP_Agmon}
S.~Agmon.
\newblock The relaxation method for linear inequalities.
\newblock {\em Canadian Journal of Mathematics}, 6(3):382--392, 1954.

\bibitem{BauschkeFixPointsNonExpansive}
H.~H. Bauschke.
\newblock The approximation of fixed points of compositions of nonexpansive
  mappings in hilbert space.
\newblock {\em Journal of Mathematical Analysis and Applications}, 202(1):150
  -- 159, 1996.

\bibitem{BauschkeOnProjection}
H.~H. Bauschke and J.~M. Borwein.
\newblock On projection algorithms for solving convex feasibility problems.
\newblock {\em SIAM Review}, 38(3):367--426, 1996.

\bibitem{bauschkeCVXanal}
H.~H. Bauschke and P.~L. Combettes.
\newblock {\em {Convex Analysis and Monotone Operator Theory in Hilbert
  Spaces}}.
\newblock Springer, 2011.

\bibitem{Boyd2004}
S.~Boyd and L.~Vandenberghe.
\newblock {\em Convex Optimization}.
\newblock Cambridge University Press, New York, NY, 2004.

\bibitem{BREGMAN1967200}
L.M. Bregman.
\newblock The relaxation method of finding the common point of convex sets and
  its application to the solution of problems in convex programming.
\newblock {\em USSR Computational Mathematics and Mathematical Physics},
  7(3):200 -- 217, 1967.

\bibitem{Combettes201555}
P.~L. Combettes and I.~Yamada.
\newblock Compositions and convex combinations of averaged nonexpansive
  operators.
\newblock {\em Journal of Mathematical Analysis and Applications},
  425(1):55--70, 2015.

\bibitem{cvxpy}
S.~Diamond and S.~Boyd.
\newblock {CVXPY}: A {P}ython-embedded modeling language for convex
  optimization.
\newblock {\em Journal of Machine Learning Research}, 2016.
\newblock To appear.

\bibitem{DouglasRachford}
J.~Douglas and H.~H. Rachford.
\newblock On the numerical solution of heat conduction problems in two and
  three space variables.
\newblock {\em Trans. Amer. Math. Soc.}, 82:421--439, 1956.

\bibitem{gisSIAM2015}
P.~Giselsson.
\newblock Tight global linear convergence rate bounds for {D}ouglas-{R}achford
  splitting.
\newblock 2015.
\newblock Submitted. Available: {\url{http://arxiv.org/abs/1506.01556}}.

\bibitem{gis_line_search}
P.~Giselsson, M.~F{\"a}lt, and S.~Boyd.
\newblock Line search for averaged operator iteration.
\newblock Available: \url{http://arxiv.org/abs/1603.06772v2}, 2016.

\bibitem{cvx_v3}
M.~Grant and S.~Boyd.
\newblock {CVX}: Matlab software for disciplined convex programming, version
  3.0.
\newblock \url{http://cvxr.com/cvx}, 2016.

\bibitem{GAP_Gubin}
L.~G. Gubin, B.~T. Polyak, and E.~V. Raik.
\newblock The method of projections for finding the common point of convex
  sets.
\newblock {\em USSR Computational Mathematics and Mathematical Physics},
  7(6):1--24, 1967.

\bibitem{LionsMercier1979}
P.~L. Lions and B.~Mercier.
\newblock Splitting algorithms for the sum of two nonlinear operators.
\newblock {\em SIAM Journal on Numerical Analysis}, 16(6):964--979, 1979.

\bibitem{Nocedal}
J.~Nocedal and S.~Wright.
\newblock {\em Numerical optimization}.
\newblock Springer series in operations research and financial engineering.
  Springer, New York, NY, 2nd edition, 2006.

\bibitem{SCS}
B.~O'Donoghue, E.~Chu, N.~Parikh, and S.~Boyd.
\newblock Conic optimization via operator splitting and homogeneous self-dual
  embedding.
\newblock {\em Journal of Optimization Theory and Applications}, 2016.

\bibitem{convexjl}
M.~Udell, K.~Mohan, D.~Zeng, J~Hong, S.~Diamond, and S.~Boyd.
\newblock Convex optimization in {J}ulia.
\newblock {\em SC14 Workshop on High Performance Technical Computing in Dynamic
  Languages}, 2014.

\bibitem{ye1994}
Y.~Ye, M.~J. Todd, and S.~Mizuno.
\newblock An o($\sqrt{n}l$)-iteration homogeneous and self-dual linear
  programming algorithm.
\newblock {\em Mathematics of Operations Research}, 19(1):53--67, 1994.

\end{thebibliography}

\begin{appendix}
\section{Appendix}
\subsection{Proof of Proposition~\ref{prop:prop1}}\label{app:prop1}

We start with the first claim.
We know from~\cite[Proposition 4.8]{bauschkeCVXanal} that $\Pi_{C_i}$ is firmly nonexpansive,
and since $\alpha_i\in(0,2)$ we know from~\cite[Corollary 4.29]{bauschkeCVXanal}
that $P_{C_i}^{\alpha_i}$ are $\tfrac{\alpha_i}{2}$-averaged.

The composition $P_{C_p}^{\alpha_p}\dots P_{C_1}^{\alpha_1}$ is therefore $\beta$-averaged with $\beta$ in~\eqref{eq:beta},
according to~\cite{Combettes201555, gisSIAM2015}.
Therefore
\begin{align*}
T&=(1-\alpha)\id+\alpha((1-\beta)\id+\beta S)\\
&=(1-\alpha\beta)\id+\alpha\beta S.
\end{align*}
where $S$ is nonexpansive. Since $\alpha\in(0,\tfrac{1}{\beta})$ we have $\alpha\beta\in(0,1)$ and the first claim is proven.

The second claim is proven by noting that $P_{C_i}^{\alpha_i}$
is nonexpansive when $\alpha_i=2$~\cite[Corollary 4.10]{bauschkeCVXanal}.
This implies that the composition is nonexpansive and the claim follows directly since $\alpha\in(0,1)$.

\subsection{Proof of Proposition~\ref{prp:fixT}}\label{app:fixT}
To show this, we need the following lemma.
\begin{lem}
Suppose that $C$ is a nonempty closed and convex set and $\alpha\neq 0$.
Then $\fix P^\alpha_C=C$, i.e. $x\in C$ if and only if $P_{C}^{\alpha}x=x$.
\label{lem:rel_proj_equiv}
\end{lem}
\begin{pf}
It holds for projection operators with $\alpha=1$~\cite[Equation 4.8]{bauschkeCVXanal}.
Since
\begin{align*}
P_C^{\alpha}x=\alpha \Pi_Cx+(1-\alpha)x=x+\alpha(\Pi_Cx-x)
\end{align*}
we have $P_C^{\alpha}x=x$ if and only if $\Pi_Cx=x$ if $\alpha\neq 0$. This concludes the proof of the lemma.
\end{pf}
The result follows directly for the case A1 from~\cite[Corollary 4.37]{bauschkeCVXanal}
since $\fix P^{\alpha_i}_{C_i}=C_i$ and $P^{\alpha_i}_{C_i}$ are $\alpha_i$-averaged operators.

For the case A2, let $j$ be the index with $\alpha_j=2$ and first assume that $j\neq 1,j\neq p$.
Define
$S_1 = P_{C_p}\dots P_{C_{j+1}}$ and $S_2 = P_{C_{j-1}}\dots P_{C_{1}}$.

Since all $P_{C_i}^{\alpha_i}$ are averaged for $i=j+1,\dots,p$,
and since $\fix P_{C_i}^{\alpha_i}=C_i$ from Lemma~\ref{lem:rel_proj_equiv},
~\cite[Corollary 4.37]{bauschkeCVXanal} gives that
$S_1$ is strictly quasi-nonexpansive and that
$\fix S_1=
\cap_{i=j+1}^{p} C_i$.
The same argument shows that $S_2$
is stricty quasi-nonexpansive with $\fix S_2=\cap_{i=1}^{j-1} C_i$.

Let $T_1 = S_1P_{C_j}^2$.
Nonexpansiveness of $P_{C_j}^2$ implies quasi-nonexpansiveness,
so $T_1$ is also quasi-nonexpansive with $\fix T_1 = \cap_{i=j}^{p} C_i$ by~\cite[Proposition~4.35]{bauschkeCVXanal}.
Again applying \cite[Proposition~4.35]{bauschkeCVXanal} to $T=T_1S_2$ gives the result.

In the special case where $j=p$ or $j=1$,
the results follows in the same way for $T=P_{C_p}^2S_2$ or $T=S_1P_{C_1}^2$ respectively.

\subsection{Proof of Theorem~\ref{thm:proj-cheap}}
For the new iterate $x^{k+1}$ in \eqref{eq:xk+1Alt2}, we have $x^{k+1}=\Pi_C(x^k+\alphak r^k)$,
and therefore $x^{k+1}\in C$. Let the shortest distance to a set $\Omega$ be denoted $\rm{dist}_\Omega(x):={\|\Pi_\Omega x-x\|_2}$.
The norm of the residual then simplifies to
\begin{align}
\|r(x^{k+1})\| &= \|P^{\alpha_2}_DP^{\alpha_1}_C x^{k+1}-x^{k+1}\| \\
& = \|P^{\alpha_2}_D x^{k+1}-x^{k+1}\| \\
& = \|\alpha_2\Pi_D x^{k+1}+(1-\alpha_2)x^{k+1}-x^{k+1}\| \\
\label{eq:cheapSearch}& = \alpha_2\|\Pi_Dx^{k+1}-x^{k+1}\| \\
& = \alpha_2{\rm{dist}}_D(x^{k+1}) \\
& = \alpha_2{\rm{dist}}_D\left(\Pi_C(x^k+\alphak r^k)\right).
\end{align}
Since $C$ is affine, so is $\Pi_C$. This implies that
$\Pi_C(x^k+\alphak r^k)$ is affine in $\alphak$.
So the norm of the residual is the composition between the convex function ${\rm{dist}}_D$
and an affine function in $\alphak$, hence convex~\cite{Boyd2004} in $\alphak$.

\end{appendix}
\end{document}